\documentclass[12pt]{article}
\usepackage{pstricks,subfigure,graphicx}
\usepackage{tikz}
\usetikzlibrary{snakes,arrows}
\usepackage{amsmath,amssymb}
\usepackage{tikz}
\usepackage{color}
\usepackage{float}
\newtheorem{theorem}{Theorem}[section]
\newtheorem{lemma}[theorem]{Lemma}

\newtheorem{proposition}[theorem]{Proposition}

\newcommand{\proof}{\noindent{\bf Proof.\ }}
\newcommand{\qed}{\hfill $\square$ \bigskip}
\newcommand{\GG}{G\overline{G}}
\newcommand{\TT}{T\overline{T}}

\DeclareMathOperator {\con} {con}

\DeclareMathOperator {\diam} {diam}
\DeclareMathOperator {\rad} {rad}

\let\deg\relax
\DeclareMathOperator {\deg} {deg}

\begin{document}

\title{On the convexity number of complementary prisms of trees}

\author{
	Neethu P. K.$^{a}$
	\and
	Ullas Chandran S. V. $^{a}$
	}

\date{\today}

\maketitle

\begin{center}
	 Department of Mathematics, Mahatma Gandhi College, Kesavadasapuram,  Thiruvananthapuram-695004, Kerala, India \\
	{\tt p.kneethu.pk@gmail.com}, {\tt svuc.math@gmail.com} 
	\medskip

\end{center}
\begin{abstract}
A set of vertices $S$ of a graph $G$ is a (geodesic) convex set, if $S$ contains all the vertices belonging to any shortest path connecting between two vertices of $S$. The cardinality of maximum proper convex set of $G$ is called the convexity number, $\con(G)$, of $G$. The complementary prism $\GG$ of $G$ is obtained from the disjoint union of $G$ and its complement $\overline{G}$ by adding the edges of a perfect matching between them. In this work, we examine the convex sets of the complementary prism of a tree and derive formulas for the convexity numbers of the complementary prisms of all trees. 
\end{abstract}

\medskip\noindent

{\bf AMS Subj.\ Class:} 05C12; 05C69.

\section{Introduction}
\label{sec:intro}

A family ${\cal C}$ of subsets of a finite set $X$ is a \emph{convexity} on $X$ if $\emptyset \in {\cal C}$, $X\in {\cal C}$ and ${\cal C}$ is closed under arbitrary intersections. A subset $T$ of $X$ is a $\mathcal{C}$-\emph{convex set} if  $ T \in \mathcal{C}$. An extensive survey of abstract convexity and related combinatorial geometry can be found in ~\cite{edelman-1985,farber-1986,van-1993}. Different convexities associated with the vertex set of a graph are well-known. The most natural convexities in graphs are path convexities defined by a family of paths $\mathcal{P}$, in a way that a set $T$ of vertices of $G$ is \emph{convex} if and only if each vertex that lies on an $(u,v)$-path of $\mathcal{P}$ belongs to $T$. An extensive survey of different types of path convexities can be found in \cite{pelayo-2015}.
 
 In this paper, we consider the geodetic convexity in graphs. In this convexity, $\cal P$ is the family of geodesics (shortest paths) of the graph. A set of vertices $S$ of a graph $G$ is a  ($geodesic$) \emph{convex set}, if $S$ contains all the vertices belonging to any shortest path between two vertices of $S$. The cardinality of a maximum proper convex set of $G$ is called the \emph{convexity number}, $\con(G)$, of $G$.

Let us briefly recall the progress on the convexity number so far. The convexity number of a graph was introduced by G. Chartrand, C. E. Wall and P. Zhang in \cite{chartrand-2002}. In \cite{gimbel-2003}, it is proved that the decision problem associated with the convexity number is NP-complete. The convex sets and the convexity number of a graph have been further investigated in a sequence of papers \cite{ edelman-1985,farber-1986}. The convexity numbers of join, Cartesian products and lexicographic products have been further studied in \cite{ anand-2012,canoy-2002}. 

 If $G$ is a graph and $\overline{G}$ its complement, then the {\em complementary prism} $G\overline{G}$ of $G$ is the graph formed from the disjoint union of $G$ and $\overline{G}$ by adding the edges of a perfect matching between the corresponding vertices of $G$ and $\overline{G}$~\cite{haynes-2007}. For example,  $C_5\overline{C}_5$ is the Petersen graph. Solely from this particular reason, but also from many additional ones, complementary prisms were studied from different perspectives. Since the Petersen graph is a key example in the theory of edge colorings, it is no surprise that the chromatic index of complementary prisms was studied in~\cite{zatesko-2019}. Other topics studied on complementary prisms include domination~\cite{haynes-2009}, cycle structure~\cite{meieling-2015}, complexity properties~\cite{duarte-2017}, spectral properties~\cite{cardoso-2018}, hull number ~\cite{coelho-2018}, b-chromatic number~\cite{bandeli-2019}  and general position number~\cite{neethu-2020}. The convexity number of complementary prisms have been investigated in {\rm\cite{castonguay-2019}}. In {\rm\cite{castonguay-2019}}, $\con(G\overline{G})$ is determined when $G$ or $\overline{G}$ is disconnected and it is proved that the decision problem related associated to the convexity number is NP-complete even restricted to complementary prisms. In the same paper, a lower bound has been obtained when $\diam(G)\neq 3$. In this paper, we continue the investigation on the convexity number of the complementary prisms $\GG$, when $G$ is a tree. We derive formulas for the convexity number of the complementary prisms of all trees. In Section 2, we fix the notations, terminologies and discuss some preliminary results of convexity number already available in the literature. For each vertex $x$ in a tree $T$, let $e(x)$ be the number of pendent neighbors of $x$. If $G=(V(G),E(G))$ is a graph, then $n(G)=|V(G)|$. Fix $\Delta (T\backslash\{x\})=\max \{deg_T(y)$ | $y\in V(T)\backslash \{x\}\}$. Also, in a tree $T$ with diameter 4, denote by $c_t$ the unique central vertex of $T$, where central vertices of a graph are the vertices with minimum eccentricity. With these notations in hand, in Section 3, we prove the following theorem.

\begin{theorem}\label{mainthm} Let $T$ be any tree with al least three vertices. Then
	$$\con(T\overline{T})=
	\begin{cases}
	max\{n(T),2\Delta(T)+1\}$  $ \hspace{0.5cm}; if \diam(T)\geq 5 \,, \\
\hspace{1cm}	n(T)+\Delta(T)-1 \hspace{0.9cm} ;if \diam(T)=4 $ and $ deg_T(c_t)<\Delta(T) \,, \\
	\max\{n(T)+\Delta(T\backslash \{c_t\})-1,\\\hspace{0.9cm} 2\Delta(T)+1, \\ \hspace{0.9cm} n(T)+2e(c_t)-\Delta(T)+1\};  $ $ if \diam(T)=4 $ and $  deg_T(c_t)=\Delta(T)  \,, \\
\hspace{1cm}	n(T)+\Delta(T)-1 \hspace{1.1cm}; if \diam(T)=3 \,, \\
\hspace{1cm}	2n(T)-1 \hspace{2.3cm}; if \diam(T)=2.	
	\end{cases}$$
	\end{theorem}
		
\section{Preliminaries}
\label{sec:preliminary}

Graphs in this paper are finite, simple and undirected. For basic graph terminologies, we follow \cite{chartrand-2006}. Let $G = (V(G), E(G))$ be a graph. The maximum degree among all vertices of $G$ is denoted by $\Delta(G)$.  The \emph {open neighborhood} $N_G(v)$
is the set of neighbors of $v$, while the \emph {closed neighborhood} $N_G[v]$ is the open neighborhood together with the vertex $v$ itself. The \emph {closed neighborhood of a set} $S$ of vertices is $N_G[S] = \bigcup _{v\in S} N_G[v]$. If $N_{G}[v]$ is a clique then $v$ becomes an extreme vertex of $G$. The \emph {open neighborhood of a set} $S$ of vertices is $N_G(S) = (N_G[S])\backslash S$.  The {\em distance} $d_G(u,v)$ between vertices $u$ and $v$ is the length of a shortest $u,v$-path. An $u,v$-path of minimum length is also called an $u,v$-{\it geodesic}. The {\it eccentricity} of $u$ is ${\rm ecc}_G(u) = \max \{d_G(u,v) $ | $ v\in V(G)\}$. The {\em radius} and the {\em diameter} of $G$ are ${\rm rad}(G) = \min \{{\rm ecc}_G(v)$ | $v\in V(G)\}$ and $\diam(G) = \max \{ecc_G(v)$ | $v\in V(G)\}$, respectively. A star graph is a tree with diameter 2. A vertex $v$ is a \emph{central vertex} of $G$ if ${\rm ecc}_G(v)=\rad(G)$. The set of all central vertices is denoted by $C(G)$. We may simplify the above notation by omitting the index $G$ whenever $G$ is clear from the context. On the other hand, when we want to emphasize that a vertex is central in a graph $G$, we will say that it is {\em $G$-central}. The {\em interval} $I_G[u,v]$ between $u$ and $v$  is the set of vertices that lie on some $u,v$-geodesic of $G$. For $S\subseteq V(G)$ we set $I_G[S]=\bigcup_{_{u,v\in S}}I_G[u,v]$. Thus a set $S$ is convex if $I_G[S]=S$. The \emph{convex hull} $[S]_G$ of a set $S$ in $G$ is the smallest convex set containing $S$ in $G$. A set $S$ is a \emph{hull set} if $[S]_G=V(G)$.  This definition allows us to extend and study several problems from classical convexity to a finite and discrete form.

In this paper, we will also make use of the following results. 
\begin{proposition}
{\rm\cite[Proposition 3.2]{castonguay-2019}}
\label{proposition:2.1}
Let $G$ be a graph, $S\subseteq V(\GG)$, and $v_{1},v_{2},..v_{k}$ be a path in $G$, for $k\geq 2$. If $\{v_{1},\overline{v_{2}},... \overline{v_{k}}\}\subseteq [S]_{\GG}$, then $v_{k}\in [S]_{\GG}$.
\end{proposition}

\begin{theorem}
{\rm\cite[Theorem 3.9]{castonguay-2019}}
\label{thm:2.2}
Let $G$ be a graph with $\diam(G)\neq 3$, then $\con(\GG)\geq n(G)$.
\end{theorem}

{\begin{theorem}
{\rm\cite[Theorem 3.3]{castonguay-2019}}
\label{thm:2.3}
Let $G$ be a disconnected graph and $k$ be the order of a minimum component of $G$. Then $\con(\GG)= 2n(G)-k$.
\end{theorem}
}

\section{Proof of Theorem \ref{mainthm}}
\label{sec:Trees}
To prove the announced theorem, some preparations are needed.  
\begin{lemma}
\label{lemma:3.1}
Let $G$ be a connected graph. Then  $V(\overline{G})\cup \{u\}$ is a hull set of $V(\GG)$ for all $u \in V(G)$.
\end{lemma}
\proof
Let $S= V(\overline{G})\cup \{u\}$ and let $v$ be any vertex in $G$ distinct from $u$. Let $P:u=v_{1},v_{2},...,v_{k}=v$ be some $u,v$-path in $G$. Since $v_{1},\overline{v_{2}},\overline{v_{3}},...,\overline{v_{k}}\in [S]_{\GG}$, by Proposition~\ref{proposition:2.1}, $v\in [S]_{\GG}$. Hence $[S]_{\GG}=V(\GG)$.
\qed
\begin{lemma}
\label{lemma:3.2}
Let $T$ be a tree with $\diam(T)\geq 3$. If $u,v \in V(T)$ with $d_{T}(u,v)=3$, then $[\{u,v\}]_{\TT}=V(\TT)$.
\end{lemma}
\proof
 Let $S=\{u,v\}$ and let $P:u,x,y,v$ be a $u,v$-geodesic of length 3 in $T$. Then both $P$ and the path $Q:u, \overline{u}, \overline{v},v$ are $u,v$-geodesics in ${\TT}$. This shows that $u,x,y,v,\overline{u},\overline{v}\in I_{\TT}[S]$. Also, since the paths $P_{1}:x,\overline{x},\overline{v}$ and $P_{2}:y,\overline{y},\overline{u}$ are geodesics in $\TT$, we have that $\overline{x},\overline{y}\in I^{2}_{\TT}[S]$. Now, let $z$ be any vertex in $T$ such that $z\notin V(P)$. Since $T$ is a tree, it follows that $z$ must be non-adjacent to at least two adjacent vertices in $P$, say $u$ and $x$. This shows that $\overline{z}\in I_{\TT}[\overline{u},\overline{x}]\subseteq [S]_{\TT}$. Thus $V(\overline{T})\subseteq [S]_{\TT}$.
 Now, it follows from Lemma~\ref{lemma:3.1} that $V(\TT)\subseteq [V(\overline{T})\cup \{v\}]_{\TT}\subseteq [S]_{\TT}$.
 \qed

\begin{lemma}
\label{lemma:3.3} Let $T$ be a tree of $\diam(T)\geq 6$. Then $V(\overline{T})\subseteq [\overline{u},\overline{v}]_{\TT}$ for any pair of adjacent vertices $u,v$ in $V(T)$. 
\end{lemma}
\proof
Let $u, v$ be adjacent vertices in $T$ and let $H=V(T)\backslash ({N_{T}(u)}\cup {N_{T}(v)})$. Then it is clear that $\overline{H}\subseteq I_{{\TT}}[\overline{u},\overline{v}]$, where $\overline{H}$ denotes the corresponding vertices of $H$ in $\overline{T}$. On the other hand, choose $x$ in $N_{T}(u)$ or $y\in N_{T}(v)$, since $\diam(T)\geq 6$, it follows that there exist adjacent vertices $w,z\in H$ such that either $xw,xz\notin E(T)$ or $yw,yz\notin E(T)$, say $x$ is non-adjacent to $w$ and $z$ in $T$. This shows that $\overline{x}\in I_{{\TT}}[\overline{w},\overline{z}]\subseteq [H]_{\TT}\subseteq [\overline{u},\overline{v}]_{\TT}$. This is in turn implies that $\overline{N_{T}(u)}\subseteq [\overline{u},\overline{v}]_{\TT}$. Also, for each $x\in N_{T}(u)\backslash \{v\}$, $\overline{N_{T}(v)}\subseteq I_{{\TT}}[\overline{x},\overline{u}]\subseteq [H]_{\TT}\subseteq [\overline{u},\overline{v}]_{\TT}$. Thus we have proved that $V(\overline{T})\subseteq [\overline{u},\overline{v}]_{\TT}$.
\qed
\begin{lemma}
\label{lemma:3.4}
Let $T$ be a tree with $\diam(T) = 5$. If $u,v$ are two adjacent vertices of $T$ such that $u\in {\rm Per}(T)$, then $V(\overline{T})\subseteq [\overline{u},\overline{v}]_{\TT}$.
\end{lemma}
\proof
Let $u,v$ be two adjacent vertices of $T$ such that $u\in {\rm Per}(T)$. Let $w\in V(T)$ be such that $d_{T}(u,w)=5$. This is in turn implies that $v$ lies on the $u, w$-path of $T$. Let $P: u=u_{0},u_{1}=v,u_{2},u_{3},u_{4},u_{5}=w$ be the path from $u$ to $w$. Since $v$ is the unique neighbor of $u$ in $T$, we have that $V(\overline{T})\backslash (\overline{N_{T}(v)})\subseteq I_{{\TT}}[\overline{u},\overline{v}]$. Also $\overline{N_{T}(v)}\subseteq I_{\TT}[\overline{u_{4}},\overline{w}]\subseteq [\overline{u},\overline{v}]_{\TT}$. Thus $V(\overline{T})\subseteq [\overline{u},\overline{v}]_{\TT}$.
\qed
\begin{lemma}
\label{lemma:3.5}
Let $T$ be a tree with $\diam(T) = 5$. If $u,v,w$ be three vertices of $T$ such that at least one of them is adjacent to remaining two vertices. Then $V(\overline{T})\subseteq [\{\overline{u},\overline{v},\overline{w}\}]_{\TT}$.
\end{lemma}
\proof
Let $S=\{\overline{u},\overline{v},\overline{w}\}$ and $u,v,w\in V(T)$ be such that $uv,vw\in E(T)$. Since $T$ is a tree, we can easily verify that $V(\overline{T})\backslash (\overline{N_{T}(u)}\cup \overline{N_{T}(v)}\cup \overline{N_{T}(w)})\subseteq I_{T\overline{T}}[\overline{u},\overline{v}]\subseteq [S]_{\TT}$ and  $\overline{N_{T}(u)}\subseteq I_{T\overline{T}}[\overline{v},\overline{w}]\subseteq [S]_{\TT}$. Similarly, $\overline{N_{T}(w)}\subseteq I_{T\overline{T}}[\overline{u},\overline{v}]\subseteq [S]_{\TT}$. Now we will show that $\overline{N_{T}(v)}\subseteq [S]_{\TT}$. Let $z\in N_{T}(v)$. If $z\notin C(T)$, then ${\rm ecc}_{T}(z)\geq 4$ and thus there exists adjacent vertices $x$ and $y$ in $V(T)\backslash N_{T}(v)$ such that $\overline{z}\in I_{\TT}[\overline{x},\overline{y}]\subseteq[S]_{\TT}$. If $z\in C(T)$, then since $|C(T)|\leq 2$ and $C(T)$ is a clique, it is clear that $z$ is the unique central vertex in $N_T(v)$. But we have that $V(\overline{T})\backslash \{\overline{z}\}\subseteq [S]_{\TT}$ and ${\rm ecc}_{T}(z)=3$. Thus using a parallel argument we infer that $\overline{z}\in I_{\TT}[\overline{x},\overline{y}]\subseteq [S]_{\TT}$ for some  $x,y$ in $V(T)$. 
\qed

{\begin{lemma}
\label{lemma:3.6}
Let $T$ be a tree. Then $\con(T\overline{T})\geq 2\Delta(T)+1$.
\end{lemma}
\proof
Let $x\in V(T)$ be such that $\deg_{T}(x)=\Delta(T)$. Consider the set $H=N_{T}[x] \cup \overline{N_{T}(x)}$ (as indicated in bold in \figurename{ \ref{fig:3.1}}). We claim that $H$ is convex in $\TT$. Since $N_{T}[x]$ induces a star in $T$; and  $\overline{N_{T}(x)}$ induces a clique in $\overline{T}$, it is clear that both $N_{T}[x]$, $\overline{N_{T}(x)}$ are convex sets of $\TT$. Now, for any two distinct vertices $u,v\in N_{T}(x)$, $I_{\TT}[u,\overline{v}]=\{u,\overline{u},\overline{v}\}\subseteq H$. Also, for any $u\in N_{T}(x)$, $I_{\TT}[x,\overline{u}]=\{x,u,\overline{u}\}\subseteq H$. Thus $H$ is convex in $\TT$ and hence $\con(T\overline{T})\geq |H|=2\Delta(T)+1$}. 
\qed
\renewcommand{\thefigure}{3.1}
\begin{figure}[H]
\begin{center}
\begin{tikzpicture}[scale=1.6]
\tikzstyle{every node}=[draw,shape=circle,minimum size=16mm,font=\fontsize{7}{7}\sffamily];
\node[ultra thick,dashed](A) at(0,-1){$V(T)\backslash N_{T}(x)$};
\node[ultra thick](B) at(0,1){$N_{T}(x)$};
\node[ultra thick,dashed](A1) at(2.5,-1){$\overline{V(T)\backslash N_{T}(x)}$};
\node[ultra thick](B1) at(2,1){$\overline{N_{T}(x)}$};
\node[ultra thick,minimum size=0.1mm](x) at(0,0){${x}$};
\node[ultra thick,minimum size=0.1mm,dashed](x1) at(1.5,0){$\overline{x}$};
\node[draw=none,font=\small] (T) at (0,-2.1){$T$};
\node[draw=none,font=\small] (H) at (2.5,-2.1){$\overline{T}$};
\path[-][dashed](A) edge [bend left] (B);
\draw (B)-- (x);
\draw[dashed] (A)-- (A1);
\draw[dashed] (B)-- (B1);
\draw (x)-- (x1);
\draw (A1)-- (x1);
\draw[dashed]  (A1)--(B1);
\draw(-.75,1.75)rectangle(.75,-1.75);
\draw(1.1,1.75)rectangle(3.2,-1.75);
\end{tikzpicture}
\end{center}
\caption{Illustration of the set $H$ in Lemma \ref{lemma:3.6}}
\label{fig:3.1}
\end{figure}

\begin{theorem}
\label{thm:3.7}
If $T$ is a tree with  $\diam(T)\geq 5$, then 
\begin{center} $\con(T\overline{T})= max\{n(T),2\Delta(T)+1\}$.
	\end{center}
\end{theorem}
\textbf{Proof}:
 Theorem~\ref{thm:2.2} and Lemma~\ref{lemma:3.6} implies that $\con(T\overline{T})\geq max\{n(T),2\Delta(T)+1\}$. Hence in the following we prove that $\con(T\overline{T})\leq max\{n(T),2\Delta(T)+1\}$. For, let $S$ be any proper convex set in $T\overline{T}$. Fix $M=S\cap V(T)$ and $\overline{N}=S\cap V(\overline{T})$. If $M=\emptyset$ or $\overline{N}=\emptyset$, then there is nothing to prove. So, assume that $M\neq \emptyset$ and $\overline{N}\neq \emptyset$. Now, since $M\neq \emptyset$, it follows from Lemma~\ref{lemma:3.1} that $|\overline{N}|\leq n(T)-1$. If $|M|=1$, then $|S|\leq n(T)-1$. Hence assume that $|M|\geq 2$. If there exists $u,v\in M$ with $d_{T}(u,v)=3$ then by Lemma~\ref{lemma:3.2}, $V(\TT)=[\overline{u},\overline{v}]_{\TT}\subseteq [M]_{\TT}\subseteq S$. Hence for all $u, v\in M$ either $d_{T}(u,v)>3$ or $d_{T}(u,v)\leq2$. Now, we consider the following two cases. \\ 
\textbf{Case 1}: There exists $u,v \in M $ with $d_{T}(u,v)>3$. \\ In this case, we prove that $S$ has at most $n(T)$ vertices. We first prove the following claim.\\ 
\textbf{Claim 1:} $d_{T}(x,y)>3$ for all $x,y\in M$. \\ 
\textbf{Subclaim 1}: Let $x,y\in M$ be such that $d_{T}(x,y)>3$. Then $d_{T}(x,z)>3$ , for all $z$ in $M$ distinct from $x$.\\ 
\textbf{Subcase 1.1}: $\diam(T)=5$. \\
Assume the contrary that $d_{T}(x,y)>3 $ and $d_{T}(x,z)\leq 2$ for some  $x, y,z\in M$ with $z\neq x$. First consider the case that $d_{T}(x,z)=1$. Now, since $S$ is a proper convex set in  $\TT$, it follows from Lemma~\ref{lemma:3.2} that $d_{T}(z,y)>3$. Recall that $d_{T}(x,z)=1,$ $d_{T}(x,y)>3$ and $d_{T}(z,y)>3$. Hence either $x\in Per(T)$ or $z\in Per(T)$. Thus by Lemma~\ref{lemma:3.4}, $V(\overline{T})\subseteq [S]_{\TT}=S.$ Thus it follows from Lemma~\ref{lemma:3.1} that $S=V(\TT)$, a contradiction. 

Next, consider the case $d_{T}(x,z)=2$ and let $P:x,w,z$ be an $x,z$-geodesic in $T$. Then $d_T(x, w)=1$ and $w\in I_{T}[x,z]\subseteq I_{\TT}[x,z]\subseteq [S]_{\TT}$. Hence using a parallel argument we infer that $S=V(\TT)$, a contradiction.\\  
\textbf{Subcase 1.2}: $\diam(T)\geq 6$. \\
Assume the contrary that $d_{T}(x,y)>3 $ and $d_{T}(x,z)\leq 2$ for some  $x, y,z\in M$ with $z\neq x$. First consider the case that $d_{T}(x,z)=1$. Then  $d_{T}(z,y)\geq 3$. This shows that the paths on vertices $x,\overline{x},\overline{y},y$ and $z,\overline{z},\overline{y},y$ are geodesics in $\TT$. Hence the vertices $\overline{x}$ and $\overline{z}$ belong to $\overline{N}$ and so by Lemma~\ref{lemma:3.3}, $V(\overline{T})\subseteq [\overline{x},\overline{z}]_{\TT}\subseteq [S]_{\TT}=S.$ Hence it follows from Lemma~\ref{lemma:3.1} that $S=V(\TT)$, a contradiction.

Next, consider the case $d_{T}(x,z)=2$ and let $P:x,w,z$ be an $x,z$-geodesic in $T$. Then $d_T(x,w)=1$ and $w\in I_{T}[x,z]\subseteq I_{\TT}[x,z]\subseteq [S]_{\TT}=S$. Hence using the above parallel arguments we infer that $S=V(\TT)$, a contradiction. Hence in all cases, the Subclaim 1 follows. Next, we deduce Claim 1.

Let $x$ and $y$ be two vertices in $M$. Then it follows from Lemma~\ref{lemma:3.2} that $d_{T}(x,y)\neq 3$. Recall that $d_{T}(u,v)>3$. By Subclaim 1, $d_{T}(u,x)>3$ and $d_{T}(x,v)>3$. Again by applying Subclaim 1, we have that $d_{T}(x,y)>3$ and hence the Claim 1 follows. \\ 
\textbf{Claim 2:} $|\overline{N}|\leq n(T)-|M|$.\\
Let $N=\{x\in V(T)|$ $\overline{x}\in \overline{N}\}$. Then $|N|=|\overline{N}|$. Suppose that $N$ intersects $N_{T}(M)$, say $y\in N_{T}(M)\cap N$. Let $x\in M$ be a neighbor of $y$. Then $y\in I_{T\overline{T}}[x,\overline{y}]\subseteq S$. This shows that $y\in M$  contradicts  Claim 1.  Hence $N\cap N_{T}(M)=\emptyset$. Again by Claim 1 there is no common neighbours between any two vertices in $M$. This shows that $|N_{T}(M)|\geq |M|$. Since $N\cap N_{T}(M)=\emptyset$, we have that $|N|\leq n(T)-|N_{T}(M)|\leq n(T)-|M|$. We have thus proved Claim 2 and hence $|S|\leq n(T)$. \\ 
\textbf{Case 2:} $d_{T}(u,v)\leq 2$, for all $u,v$ in $M$.\\ In this case we prove that $|S|\leq 2\Delta(T)+1$. Since $d_T(u, v)\leq 2$, for all $u, v\in M$ and $M$ is convex in $\TT$, one can easily observe that $M$ is also a convex set in $T$. This shows that $M$ induces a star in $T$ and so $|M|\leq \triangle(T)+1$. Now, since $M$ induces a star in $T$ and $|M|\geq 2$, we can choose adjacent vertices $x_t$ and $u_t$ in $M$. In the following, we claim that $\overline{N}\subseteq \overline{M}$.\\ 
\textbf{Subcase 2.1}: $\diam(T)=5$. \\
Assume to the contrary that there exists $\overline{y}\in \overline{N}$ such that $y \notin M$. Then we have that $d_{T}(x_t,y)\geq 2$. Otherwise $y\in I_{\TT}[x_t,\overline{y}]$. Similarly, $d_{T}(u_t,y)\geq 2$. This shows that $\overline{x_t} \in I_{T\overline{T}}[\overline{y},x_t]\subseteq [S]_{\TT}$ and $\overline{u_t} \in I_{T\overline{T}}[\overline{y},u_t]\subseteq [S]_{\TT}.$ Recall that $S$ is a convex set in $\TT$ and so it is clear that $\overline{x_t},\overline{u_t} \in \overline{N}$. Thus using   similar arguments for the remaining vertices of $M$,  we can conclude that $\overline{M}\subseteq \overline{N}$. Now,  if $x_t\in Per(T)$ or $u_t\in Per(T)$, then Lemma~\ref{lemma:3.4} implies that $V(\overline{T})\subseteq [\overline{x_t},\overline{u_t}]_{\TT}\subseteq \overline{N}$. Hence  Lemma~\ref{lemma:3.1} leads to the fact that $[S]_{\TT}=S=V(\TT)$, a contradiction. On the other hand, assume that both $x_t$ and $u_t$ are not peripheral vertices. Then since $\diam(T)=5$, it is clear that $|N_{T}(u_t,x_t)|\geq 2$. Since $u_t,x_t\in N$ and $S$ is a proper convex set of $\TT$, Lemma~\ref{lemma:3.5} leads to the fact that $N_{T}(u_t,x_t)\cap N=\emptyset$. Now, recall that $M$ induces a star in $T$ and so $|M|=2$. Moreover, $|N|\leq n(T)-2$. This shows that $|S|\leq (n(T)-2)+2=n(T)$, a contradiction. We have thus proved the claim $\overline{N}\subseteq \overline{M}$. Now, let $x_t$ be the central vertex of the induced subgraph of $M$ in $T$. Again Lemmas \ref{lemma:3.1} and \ref{lemma:3.5} show that $\overline{N}\subsetneq \overline{M}$.  Hence $|\overline{N}|\leq \triangle(T)$ and so $|S|\leq |M|+|\overline{N}|\leq 2\triangle(T)+1$. \\ 
\textbf{Subcase 2.2}:  $\diam(T)\geq 6$. \\
Assume to the contrary that there exists $\overline{y}\in \overline{N}$ such that $y \notin M$. Now, if $x_t$ and $y$ are adjacent in $T$, then $y\in I_{\TT}[x_t, \overline{y}]\subseteq S$. This leads to the fact that $y\in M$, impossible. Thus we have that $d_{T}(x_t,y)\geq 2$. Similarly, $d_{T}(u_t,y)\geq 2$. This shows that $\overline{x_t} \in I_{T\overline{T}}[\overline{y},x_t]\subseteq S$ and $\overline{u_t} \in I_{T\overline{T}}[\overline{y},u_t]\subseteq S.$ Recall that $S$ is a convex set. Hence $\overline{x_t},\overline{u_t} \in \overline{N}$. Moreover, since $\overline{x_t}\overline{u_t}\notin E(\overline{T})$, it follows from Lemma~\ref{lemma:3.3} that $V(\overline{T})\subseteq [\overline{x_t},\overline{u_t}]_{\TT}=[\overline{x_t},\overline{u_t}]_{\overline{T}}\subseteq S$. Now by Lemma~\ref{lemma:3.1}, we must have $S=V(\TT)$, a contradiction. Thus $\overline{N}\subseteq \overline{M}$. Again, it follows from Lemma~\ref{lemma:3.3} that $\overline{N}$ induces a clique in $\overline{T}$. Also recall that $M$ induces a star in $T$. Thus $|\overline{N}|\leq \triangle(T)$ and so $|S|\leq |M|+|\overline{N}|\leq 2\triangle(T)+1$.
\qed 

\begin{lemma}
\label{lemma:3.8}
Let $T$ be a tree with $\diam(T) = 4$. If $u,v,x$ are three vertices of $T$ such that $x\in C(T)$ and $u,v\in N_{T}(x)$ with $\deg_{T}(u)\geq 2$. Then $V(\overline{T})\subseteq [\overline{u},\overline{v},\overline{x}]_{\TT}$.
\end{lemma}
\proof
Let $S=\{\overline{u},\overline{v},\overline{x}\}$. Since $T$ is a tree, one can easily verify that $V(\overline{T})\backslash (\overline{N_{T}(u)}\cup \overline{N_{T}(x)})\subseteq I_{T\overline{T}}[\overline{u},\overline{x}]\subseteq [S]_{\TT}$ and $\overline{N_{T}(u)}\subseteq I_{T\overline{T}}[\overline{v},\overline{x}]\subseteq [S]_{\TT}$. Also, since $\deg_{T}(u)\geq 2$, there exists at least one vertex $z\in N_{T}(u)$ distinct from $x$. Then $\overline{N_{T}(x)}\subseteq I_{T\overline{T}}[\overline{u},\overline{z}]\subseteq [S]_{\TT}$ and so  $V(\overline{T})\subseteq [\overline{u},\overline{v},\overline{x}]_{\TT}$.
\qed
{\begin{lemma}
\label{lemma:3.9}
Let $T$ be a tree of diameter at most 4 and let $w$ be a non-pendant vertex of $T$. Then $\con(T\overline{T})\geq n(T)-\deg_{T}(w)+2e(w)+1$.
\end{lemma}
\proof
Let $B$ be the set of all non-pendant neighbors of $w$ in $T$.  Consider the set $H=(N_{T}[w]\backslash B)\cup (V(\overline{T})\backslash \overline{B})$ (as indicated in bold in \figurename{ \ref{fig:3.2}}). Then $|H|=e(w)+1+n(T)-(\deg_{T}(w)-e(w))= n(T)+2e(w)-\deg_{T}(w)+1$. We claim that $H$ is a convex set of $\TT$. Since $(N_{T}[w]\backslash B)$ induces a star in $T$, it is a convex set in $\TT$. First, if $\diam(T)\leq 3$, then the only non-pendant vertices are the central vertices of $T$. Since  $\rad(T)\leq 2$, it follows that $V(\overline{T})\backslash \overline{B}$ is a convex  set in $\overline{T}$. Next, suppose that $\diam(T)=4$ and $w\notin C(T)$. Then the unique central vertex $c_t$ is the only non-pendant neigbor of $w$ in $T$. Now, since $\overline{c_t}$ is an extreme vertex of $\overline{T}$, $V(\overline{T})\backslash \overline{B}$ is a convex set $\TT$. Now, if $w$ is a central vertex, then $V(\overline{T})\backslash( \overline{B}\cup \{\overline{w}\})$ is a clique. Note that, for any $z\in N_{T}(w)\backslash B,$ $I_{\TT}[\overline{w},\overline{z}]\subseteq V(\overline{T})\backslash \overline{B}$. Thus $(V(\overline{T})\backslash \overline{B})$ is a convex set in $\overline{T}$. Since $\diam(\overline{T})=2$, it is also convex in $\TT$. Now, for $v\in N_{T}[w]\backslash B$ and $\overline{u}\in V(\overline{T})\backslash \overline{B}$, we have that $I_{\TT}[v,\overline{u}]=\{v,\overline{v},\overline{u}\}\subseteq H$. Thus we have proved that $H$ is a convex set of $\TT$. Hence $\con(\TT)\geq |H|= n(T)-\deg_{T}(w)+2e(w)+1$.}
\qed
\renewcommand{\thefigure}{3.2}
\begin{figure}[H]
\begin{center}
\begin{tikzpicture}[scale=1.6]
\tikzstyle{every node}=[draw,shape=circle,minimum size=16mm,font=\fontsize{7}{7}\sffamily];
\node[ultra thick](A) at(0,-1){$N_{T}(w)\backslash B$};
\node[ultra thick,dashed](B) at(0,1){$B$};
\node[ultra thick](A1) at(2,-1){$\overline{N_{T}(w)\backslash B}$};
\node[ultra thick,dashed](B1) at(2,1){$\overline{B}$};
\node[ultra thick,minimum size=0.1mm](x) at(0,0){${w}$};
\node[ultra thick,minimum size=0.1mm](x1) at(1.5,0){$\overline{w}$};
\node[draw=none] (ellipsis1) at (0,1.7){$\vdots$};
\node[draw=none] (ellipsis2) at (0,-1.6){$\vdots$};
\node[ultra thick](V) at(3,0){$\overline{V(T)\backslash N_{T}[w]}$};
\node[draw=none,font=\small] (T) at (0,-2.2){$T$};
\node[draw=none,font=\small] (H) at (2.5,-2.2){$\overline{T}$};
\draw (A)-- (x);
\draw (B)-- (x);
\draw[dashed] (A)-- (A1);
\draw[dashed] (B)-- (B1);
\draw (x)-- (x1);
\draw (x1)-- (V);
\draw[dashed] (B1)-- (V);
\draw (A1)-- (V);
\draw  (A1)--(B1);
\draw(-.8,1.9)rectangle(.6,-1.9);
\draw(1.1,1.9)rectangle(3.7,-1.9);
\end{tikzpicture}
\end{center}
\caption{Illustration of the set $H$ in Lemma \ref{lemma:3.8}}
\label{fig:3.2}
\end{figure}

\begin{theorem}
\label{thm:3.10}
Let $T$ be a tree of diameter $4$ and let $c_t$ be the unique T-central vertex. Then
$$\con(T\overline{T})=
\begin{cases}
 n(T)+\Delta(T)-1;  $ $ $ if $ deg_T(c_t)<\Delta(T) \,, \\
  \max\{n(T)+\Delta(T\backslash \{c_t\})-1, 2\Delta(T)+1,\\ \hspace{4cm} n(T)+2e(c_t)-\Delta(T)+1\}$ $
  \hspace{0.3cm}; otherwise \,.
 \end{cases}$$ 
\end{theorem}
\proof
We consider the following two cases.\\
\textbf{Case 1}: $\deg_{T}(c_t)<\Delta(T)$.\\
Let $x_{\delta}$ be a vertex of maximum degree in $T$. Since $\diam(T)=4$ and $x_{\delta}$ is not a central vertex of $T$, we have that $x_{\delta}$ has exactly $\Delta(T)-1$ pendant neighbors in $T$. Then by Lemma~\ref{lemma:3.9}, $\con(T\overline{T})\geq n(T)-\Delta(T)+2(\Delta(T)-1)+1=n(T)+\Delta(T)-1$.
In the following, we prove that $con(\TT)\leq n(T)+\Delta(T)-1$.
Suppose that $T$ contains a proper convex set $S$ with $|S|> n(T)+\Delta(T)-1$ in $T\overline{T}$. Fix $M=S\cap V(T)$ and $\overline{N}=S\cap V(\overline{T})$. If $M=\emptyset$ or $\overline{N}=\emptyset$, then $|S|\leq n(T)$. So, assume that $M\neq \emptyset$ and $\overline{N}\neq \emptyset$. Also since $M\neq \emptyset$, it follows from Lemma~\ref{lemma:3.1} that $|\overline{N}|\leq n(T)-1$. If $|M|=1$, then $|S|\leq n(T)-1$. Hence assume that $|M|\geq 2$. If there exists $u,v\in M$ with $d_{T}(u,v)=3$ then by Lemma~\ref{lemma:3.2}, $[M]_{\TT}=V(T\overline{T})$. Hence for $u, v\in M$, either $d_{T}(u,v)=4$ or $d_{T}(u,v)\leq 2$. We consider the following two subcases. \\ 
\textbf{Subcase 1.1}: There exists $u,v \in M $ with $d_{T}(u,v)=4$.\\ In the following, we prove  that $S$ contains at most $n(T)$ vertices.\\ 
\textbf{Claim 1:} $d_{T}(x,y)=4$, for all $x,y\in M$. \\
\textbf{Subclaim 1}: Let $x,y\in M$ be such that $d_{T}(x,y)=4$. Then $d_{T}(x,z)=4$ , for all $z$ in $M$ distinct from $x$.\\
Assume the contrary that $d_{T}(x,y)=4 $ and $d_{T}(x,z)\leq 2$ for some  $z\in M$ with $z\neq x$. First consider the case that $d_{T}(x,z)=1$. Now, since $S$ is a proper convex set in  $\TT$, it follows from Lemma~\ref{lemma:3.2} that $d_{T}(z,y)=4$. Thus we have $d_{T}(x,z)=1,$ $d_{T}(x,y)=4$ and $d_{T}(z,y)=4$, this is not possible in a tree.\\ 
Next, consider the case $d_{T}(x,z)=2$ and let $P:x,w,z$ be an $x,z$-geodesic in $T$. Then $d_T(w, x)=1$. This shows that $x, y, w\in M$ with $d_T(x, y)=4; d_{T}(x, w)=1$. Hence using the above parallel arguments, we infer a contraction. Hence Subclaim 1 follows.

Now, let $x$ and $y$ be any two vertices in $M$. Then it follows from Lemma~\ref{lemma:3.2} that $d_{T}(x,y)\neq 3$. Recall that $d_{T}(u,v)>3$. By Subclaim 1, $d_{T}(u,x)>3$ and $d_{T}(x,v)>3$. Again, by applying Subclaim 1, we have that $d_{T}(x,y)>3$ and so Claim 1 follows. \\
\textbf{Claim 2:} $|\overline{N}|\leq n(T)-|M|$.\\
Let $N=\{x\in V(T)|\overline{x}\in \overline{N}\}$. Then $|N|=|\overline{N}|$. Suppose that $N$ intersects $N_{T}(M)$, say $y\in N_{T}(M)\cap N$. Let $x\in M$ be a neighbor of $y$ in $T$. Then $y\in I_{T\overline{T}}[x,\overline{y}]$. This shows that $y\in S$ and so $y\in M$. But from Claim 1, we have that $d_{T}(u,v)>3$ , for all $u,v$ in $M$. This is a contradiction. Hence $N\cap N_{T}(M)=\emptyset$. Again by Claim 1 there is no common neighbors between any two vertices in $M$. This shows that $|N_{T}(M)|\geq |M|$. Since $N\cap N_{T}(M)=\emptyset$, we have that $|N|\leq n(T)-|N_{T}(M)|\leq n(T)-|M|$. Thus, Claim 2 follows.\\
From Claim 2,  one can easily conclude that $|S|\leq n(T)$. \\ 
\textbf{Subcase 1.2:} $d_{T}(u,v)\leq 2$, for all $u,v$ in $M$. \\ In this case we prove that $|S|\leq n(T)+\triangle(T)-1$. Since $d_T(u, v)\leq 2$, for all $u, v\in M$ and $M$ is a convex set in $\TT$, $M$ must be a convex set in $T$. This shows that $M$ induces a star in $T$ and so $|M|\leq \triangle(T)+1$. Also we have $|\overline{N}|\leq n(T)-1$. Thus $|S|= n(T)+\Delta(T)$. This is in turn implies that $M= N_T[x_{\delta}]$,  where $d(x_{\delta})=\triangle(T)$; and $N=V(\overline{T})\backslash \{\overline{w}\}$ for some extreme vertex $\overline{w}$ of $V(\overline{T})$. But the unique extreme vertex of $\overline{T}$ is $\overline{c_t}$ and $c_t \in N(x_{\delta})$. This in turn implies that for any vertex $w\in V(T)$ with $d_{T}(w,c_t)\geq 2$, $\overline{c_t}\in I_{\TT}[c_t,\overline{w}]$. This in turn implies that $S$ cannot be a convex set in $\TT$. Hence $|S|\leq n(T)+\triangle(T)-1$. \\
\textbf{Case 2:} $\deg_{T}(c_t)=\Delta(T)$. \\ 
 Choose $u_\alpha\in V(T)$ such that $\deg_{T}(u_\alpha)=\Delta(T\backslash \{c_t\})$. Since $\diam(T)=4$ and $u_\alpha$ is not a central vertex of $T$, $u_\alpha$ has exactly $\Delta(T\backslash\{c_t\})-1$ pendant neighbors in $T$. Hence Lemma~\ref{lemma:3.9} implies that 
$\con(\TT)\geq n(T)-\Delta(T\backslash \{t\})+2(\Delta(T\backslash \{c_t\})-1)+1=n(T)+\Delta(T\backslash\{c_t\})-1$.
Also it follows from Lemmas 3.6 and \ref{lemma:3.9} that $\con(\TT)\geq 2\Delta(T)+1$ and $\con(\TT)\geq n(T)-\Delta(T)+2e(c_t)+1$. Thus it remains to prove that $\con(\TT)\leq \max\{n(T)+\Delta(T\backslash \{c_t\})-1, 2\Delta(T)+1,n(T)-\Delta(T)+2e(c_t)+1\}$.

Suppose that $T$ contains a proper convex set $S$ with $|S|> n(T)+\Delta(T\backslash \{c_t\})-1$ in $T\overline{T}$ and let $M=S\cap V(T)$ and $\overline{N}=S\cap V(\overline{T})$. If $M=\emptyset$ or $\overline{N}=\emptyset$, then $|S|\leq n(T)$. So, assume that $M\neq \emptyset$ and $\overline{N}\neq \emptyset$. Now, since $M\neq \emptyset$, it follows from Lemma~\ref{lemma:3.1} that $|\overline{N}|\leq n(T)-1$. If $|M|=1$, then $|S|\leq n(T)-1$. Hence assume that $|M|\geq 2$. If there exists $u,v\in M$ with $d_{T}(u,v)=3$ then by Lemma~\ref{lemma:3.2}, $[M]_{\TT}=V(T\overline{T})$. This leads to a contradiction to the fact that $S$ is a proper convex set of $\TT$. Hence either $d_{T}(u,v)=4$ or $d_{T}(u,v)\leq 2$, for all $u,v\in M$. We consider the following two subcases. \\ 
\textbf{Subcase 2.1}: There exists $u,v \in M $ with $d_{T}(u,v)=4$.\\ Using the same arguments in Subcase 1.1, one can easily conclude that $S$ has at most $n(T)$ vertices.\\ 
\textbf{Subcase 2.2:} $d_{T}(u,v)\leq 2$, for all $u,v$ in $M$. \\ In this case we claim that $|S|\leq \max\{2\triangle(T)+1,n(T)-\Delta(T)+2e(c_t)+1\}$. Similar to the proof of Subcase 1.2, we infer that $M$ induces a star in $T$. This shows that $M\subseteq N_{T}[z]$ for some $z\in V(T)$ and so $|M|\leq \triangle(T)+1$. First, if $z \neq c_t$, then $|M|\leq |N_{T}[z]|\leq|N_{T}[u_\alpha]|$. Also, Lemma \ref{lemma:3.1} implies that $|N|\leq n(T)-1$. Thus $|S|\leq n(T)+\Delta(T\backslash \{c_t\})-1$. Now consider the case $z=c_t$. First, If $N\subseteq M$, then $|\overline{N}|\leq \Delta(T)+1$. Lemmas \ref{lemma:3.1} and \ref{lemma:3.8} implies that  $N\subsetneq M$ and so $|\overline{N}|\leq \triangle(T)$. Hence we have that $|S|\leq |M|+|\overline{N}|\leq 2\triangle(T)+1$. On the other hand,  suppose that there exists a vertex $\overline{y}\in \overline{N}$ such that $y\notin M$. Then $d_{T}(y,w)\geq 2$, for all $w\in M$. Otherwise, $y\in I_{\TT}[w,\overline{y}]\subseteq [S]_{\TT}$, impossible. This shows that $\overline{w}\in I_{\TT}[w,\overline{y}]\subseteq [S]_{\TT}$ for all $w\in M$. Thus $M\subseteq N$. Now if $w$ is a non-pendant vertex in $N_{T}(c_t)$, then $|N_{T}(w,c_t)|\geq 2$. Suppose that $w\in M$. Then since $w,x\in N$ and $S$ is a proper convex set of $\TT$, Lemma~\ref{lemma:3.8} implies that $\overline{N_{T}(w,c_t)}\cap \overline{N}=\emptyset$. Thus $|M|\leq 2$ and $|N|\leq n(T)-2$ and hence $|S|\leq (n(T)-2)+2=n(T)$, a contradiction. Hence $M$ does not contains non-pendant vertices from $N_{T}(c_t)$. Which implies that $|M|\leq e(c_t)+1$. Also for any non-pendant vertex $w$ in $N_{T}(c_t)$, $I_{\TT}[c_t,\overline{w}]$ contains $w$ implies that $N\cap B=\emptyset$ and hence $|\overline{N}|\leq n(T)-(\Delta(T)-e(c_t))$. Thus $|S|\leq e(c_t)+1+n(T)-\Delta(T)+e(c_t)=n(T)-\Delta(T)+2e(c_t)+1$. This completes the proof.
\qed
\begin{theorem}
\label{thm:3.11}
If $T$ is a tree with  $\diam(T)=3$, then $\con(T\overline{T})= n(T)+\Delta(T)-1$.
\end{theorem}
\proof
Let $x$ and $y$ be the central vertices of $T$ with $\deg_{T}(x)=\Delta(T)$.  Since $x$ has exactly $\Delta(T)-1$ pendant neighbors, Lemma~\ref{lemma:3.9} implies that
$\con(\TT)\geq  n(T)-\Delta(T)+2(\Delta(T)-1)+1=n(T)+\Delta(T)-1.$

In the following we prove that $con(\TT)\leq  n(T)+\Delta(T)-1$.
Suppose that $\TT$ contains a proper convex set $S$ with $|S|> n(T)+\Delta(T)-1$ in $T\overline{T}$. Fix $M=S\cap V(T)$ and $\overline{N}=S\cap V(\overline{T})$. If $M=\emptyset$ or $\overline{N}=\emptyset$, then $|S|\leq n(T)$. So, assume that $M\neq \emptyset$ and $\overline{N}\neq \emptyset$. Also since $M\neq \emptyset$, it follows from Lemma~\ref{lemma:3.1} that $|\overline{N}|\leq n(T)-1$. On the other hand, for any $u\in N_{T}(x)\backslash \{y\}$ and $v\in N_{T}(y)\backslash \{x\}$, we have that $d_{T}(u,v)=3$. Thus Lemma~\ref{lemma:3.2} implies that either $M\subseteq N_{T}[x]$ or $M\subseteq N_{T}[y]$. This proves that $|M|\leq \Delta(T)+1$.\\
 This in turn implies that $|N|=n(T)-1$ and so $|M|= \Delta(T)+1$. Thus $M=N_{T}[x]$ and $\overline{N}=V(\overline{T})\backslash \{\overline{w}\}$ for some $w\in V(T)$. Note that $V(\overline{T})\backslash \{\overline{x}\}$ and $V(\overline{T})\backslash \{\overline{y}\}$ are the only maximum proper convex sets in $\overline{T}$. Hence $w=x$ or $w=y$. But $x,y\in M$ and for any $u\in N_{T}(x)$, $\overline{y}\in I_{\TT}[y,\overline{u}]$. Also for any $v\in N_{T}(y)$, $\overline{x}\in I_{\TT}[x,\overline{v}]$. Thus $S$ cannot be convex a convex set in $\TT$, impossible. Hence $\con(\TT)\leq  n(T)+\Delta(T)-1$.
  \qed \\
The following theorem is an immediate consequence of Theorem~\ref{thm:2.3}
 \begin{theorem}
\label{thm:3.12}
Let $T$ be a tree of diameter 2. Then $\con(\TT)=2n(T)-1.$
\end{theorem}

\section*{Acknowledgments}
Neethu P K acknowledges the Council of Scientific and Industrial Research(CSIR), Govt. of India for providing financial assistance in the form of Junior Research Fellowship.


\end{document}